\documentclass[a4paper]{amsart}
\usepackage{amscd}
\usepackage{amsthm}
\usepackage{graphicx}
\usepackage{amsmath}
\usepackage{amsfonts}
\usepackage{amssymb}

\begin{document}
\title[A remark on certain integral operators of fractional type]{A remark on certain integral operators of fractional type}
\author{Pablo A. Rocha}
\address{Universidad Nacional del Sur, INMABB (Conicet), Bah\'{\i}a Blanca, 8000 Buenos Aires, Argentina}
\email{pablo.rocha@uns.edu.ar}
\thanks{\textbf{Key
words and phrases}: Integral operators, Hardy spaces.}
\thanks{\textbf{2.010
Math. Subject Classification}: 42B20, 42B30.}
\thanks{Partially supported by CONICET}
\maketitle
\begin{abstract}
For $m, n \in \mathbb{N}$, $1< m \leq n$, we write $n = n_1+ ...+ n_m$ where $\{ n_1, ..., n_m \} \subset \mathbb{N}$.  Let $A_1, ..., A_m$ be $n \times n$ singular real matrices such that $$\bigoplus_{i=1}^{m} \bigcap_{1\leq j \neq i \leq m} \mathcal{N}_j = \mathbb{R}^{n},$$ where
$\mathcal{N}_j = \{ x : A_j x = 0 \}$, $dim(\mathcal{N}_j)=n-n_j$ and $A_1+ ...+ A_m$ is invertible. In this paper we study integral operators of the form
$$Tf(x)= \int_{\mathbb{R}^{n}} \, |x-A_1 y|^{-n_1 + \alpha_1} \cdot \cdot \cdot |x-A_m y|^{-n_m + \alpha_m} f(y) \, dy,$$
$n = n_1+ ...+ n_m$, $\frac{\alpha_1}{n_1} = ... = \frac{\alpha_m}{n_m}=r$, $0 < r < 1$, and the matrices $A_{i}'$s are as above. We obtain the $H^{p}(\mathbb{R}^{n})-L^{q}(\mathbb{R}^{n})$ boundedness of $T$ for $0<p<\frac{1}{r}$ and $\frac{1}{q}=\frac{1}{p} - r$.
\end{abstract}

\section{Introduction}

\qquad

For $0 \leq \alpha < n$ and $m > 1$, ($m \in \mathbb{N}$), let $T_{\alpha, m}$ be the integral operator defined by
\begin{equation}
T_{\alpha, m}f(x) = \int_{\mathbb{R}^{n}} \left| x-A_1 y\right|^{-\alpha_1} \cdot \cdot \cdot \left| x-A_m y\right|^{-\alpha_m} f(y) dy, \label{Talfa}
\end{equation}
where $\alpha_1, ..., \alpha_m$ are positive constants such that $\alpha_1 + ... + \alpha_m =n- \alpha$, and $A_1, ..., A_m$ are $n \times n$ invertible matrices such that $A_i \neq A_j$ if $i \neq j$. We observe that for the case $\alpha > 0$, $m =1$, and $A_1 = I$, $T_{\alpha, 1}$ is the Riesz potential $I_{\alpha}$. Thus for $0 < \alpha < n$ the operator $T_{\alpha, m}$ is a kind of generalization of the Riesz potential. The case $\alpha = 0$ and $m > 1$ was studied under the additional assumption that $A_i - A_j$ are invertible if $i \neq j$. The behavior of this class of operators and generalizations of them on the spaces of functions $L^{p}(\mathbb{R}^{n})$, $L^{p}(w)$, $H^{p}(\mathbb{R}^{n})$ and $H^{p}(w)$ was studied in \cite{G-U1}, \cite{G-U2}, \cite{Riveros Urciuolo}, \cite{Riv U}, \cite{R-U} and \cite{R-U2}.

If $0 < \alpha < n$ and $m > 1$, then the operator $T_{\alpha, m}$ has the same behavior that the Riesz potential on $L^{p}(\mathbb{R}^{n})$. Indeed
\[
|T_{\alpha, m}f(x)| \leq C \sum_{j=1}^{m} \int_{\mathbb{R}^{n}} | A_{j}^{-1}x - y |^{\alpha - n} |f(y)| dy = C \sum_{j=1}^{m}  I_{\alpha}(|f|)( A_{j}^{-1}x), 
\]
for all $x \in \mathbb{R}^{n}$, this pointwise inequality implies that $T_{\alpha, m}$ is a bounded operator from $L^{p}(\mathbb{R}^{n})$ into $L^{q}(\mathbb{R}^{n})$ for $1 < p < \frac{n}{\alpha}$ and $\frac{1}{q} = \frac{1}{p} - \frac{\alpha}{n}$, and it is of type weak $(1, n/n-\alpha)$.

It is well known that the Riesz potential $I_{\alpha}$ is bounded from $H^{p}(\mathbb{R}^{n})$ into $H^{q}(\mathbb{R}^{n})$ for $0 < p \leq 1$ and $\frac{1}{q} = \frac{1}{p} - \frac{\alpha}{n}$ (see \cite{weiss} and \cite{krantz}). In \cite{R-U2}, the author jointly with M. Urciuolo proved the $H^{p}(\mathbb{R}^{n})-L^{q}(\mathbb{R}
^{n})$ boundedness of the operator $T_{\alpha, m}$ and we also showed that the $H^{p}(\mathbb{R})-H^{q}(\mathbb{R})$ boundedness does not hold for $0 < p \leq \frac{1}{1+ \alpha}$, $\frac{1}{q}=\frac{1}{p} - \alpha$ and $T_{\alpha, m}$ with $0 \leq \alpha < 1$, $m=2$, $A_1 = 1$, and $A_2=-1$. This is a significant difference with respect to the case $0 < \alpha < 1$, $n=m=1$, and $A_1 = 1$.

In this note we will prove that if we consider certain singular matrices in (\ref{Talfa}), then such an operator is still bounded from $H^{p}$ into $L^{q}$. More precisely, for $m, n \in \mathbb{N}$, $1< m \leq n$, we write $n = n_1+ ...+ n_m$ where $\{ n_1, ..., n_m \} \subset \mathbb{N}$, we also consider $n \times n$ singular real matrices $A_1, ..., A_m$ such that  $$\bigoplus_{i=1}^{m} \bigcap_{1\leq j \neq i \leq m} \mathcal{N}_j = \mathbb{R}^{n},$$
where $\mathcal{N}_j = \{ x : A_j x = 0 \}$, $dim(\mathcal{N}_j)=n-n_j$, $A_1+ ...+ A_m$ is invertible. Given $0 < r <1$ and $n_1, ..., n_m$ such that $n_1+ ... + n_m =n$, let $\alpha_1, ..., \alpha_m$ be positive constants such that $\frac{\alpha_1}{n_1} = ... = \frac{\alpha_m}{n_m}=r$, for such parameters we define the integral operator $T_r$  by
\begin{equation}
T_{r}f(x)= \int_{\mathbb{R}^{n}} \, |x-A_1 y|^{-n_1 + \alpha_1} \cdot \cdot \cdot |x-A_m y|^{-n_m + \alpha_m} f(y) \, dy, \label{T}
\end{equation}
where the matrices $A_{i}'$s are as above. 

We observe that the operator defined in (\ref{T}) can be written as in (\ref{Talfa}) taking there the matrices $A_{i}$'s singular. In fact, $T_r = T_{\beta, m}$ with $\beta_{i} = n_i - \alpha_i$ for each $i = 1, 2, ..., m$ and $\beta = nr$.

\qquad

 Our main result is the following

\qquad

\textbf{Theorem 1.} \textit{Let $T_r$ be the integral operator defined in $(\ref{T})$. If $0 < r < 1$, $0 < p <\frac{1}{r}$ and $\frac{1}{q}=\frac{1}{p} - r$,
then $T_r$ can be extended to an $H^{p}(\mathbb{R}^{n}) - L^{q}(\mathbb{R}^{n})$ bounded operator.}

\qquad

In Section 2 we state two auxiliary lemmas to get the main result in Section 3. We conclude this note with an example in Section 4.

\qquad

Throughout this paper, $c$ will denote a positive constant, not necessarily the same at each occurrence. The symbol $A \lesssim B$ stands for the inequality $A \leq c B$ for some constant $c$.

\section{Preliminary results}

Let $K$ be a  kernel in $\mathbb{R}^{n} \times \mathbb{R}^{n}$, we formally define the integral operator $T_{K}$ by $T_{K}f(x) = \int_{\mathbb{R}^{n}} K(x,y) f(y) dy$.

\qquad

We start with the following

\qquad

\textbf{Lemma 1.} \textit{Let $m, n \in \mathbb{N}$, with $1< m \leq n$, and let $n_1, ..., n_m$ be natural numbers such that $n_1 + ... + n_m = n$. For each $i=1, ..., m$ let $K_i$ be kernels in $\mathbb{R}^{n_i} \times \mathbb{R}^{n_i}$ such that the operator $T_{K_i}$ is bounded from $L^{p}(\mathbb{R}^{n_i})$ into $L^{q}(\mathbb{R}^{n_i})$ with $1 < p,q < \infty$, then the operator $T_{K_1 \otimes \cdot \cdot \cdot \otimes K_m}$ is bounded from $L^{p}(\mathbb{R}^{n})$ into
$L^{q}(\mathbb{R}^{n})$.}

\qquad

\textit{Proof.} Since $\mathbb{R}^{n} = \mathbb{R}^{n_1}\times \cdot \cdot \cdot \times \mathbb{R}^{n_m}$ let $x=(x^{1}, ..., x^{m}) \in \mathbb{R}^{n_1}\times \cdot \cdot \cdot \times \mathbb{R}^{n_m}$, now the operator $T_{K_1 \otimes \cdot \cdot \cdot \otimes K_m}$ is given by
$$T_{K_1 \otimes \cdot \cdot \cdot \otimes K_m}f(x)= \int_{\mathbb{R}^{n_1}\times \cdot \cdot \cdot \times \mathbb{R}^{n_m}} \, K_1(x^{1},y^{1}) \cdot \cdot \cdot K_m(x^{m}, y^{m}) f(y^{1}, ..., y^{m}) \, dy^{1} ... dy^{m}.$$
Using that the kernels $K_i$ define bounded operators for $1 \leq i \leq m$, the lemma follows from an iterative argument and the Minkowski's inequality for integrals.$\blacksquare$

\qquad

\textbf{Lemma 2.} \textit{Let $m, n \in \mathbb{N}$, with $1< m \leq n$, and let $n_1, ..., n_m$ be natural numbers such that $n_1 + ... + n_m = n$. If $A_1, ..., A_m$ are $n \times n$ singular real matrices such that $$\bigoplus_{i=1}^{m} \bigcap_{1\leq j \neq i \leq m} \mathcal{N}_j = \mathbb{R}^{n},$$ where
$\mathcal{N}_j = \{ x : A_j x = 0 \}$, $dim(\mathcal{N}_j)=n-n_j$, and $A_1+ ...+ A_m$ is invertible, then there exist two $n \times n$ invertible matrices $B$ and $C$ such that $B^{-1}A_j C$ is the canonical projection from $\mathbb{R}^{n}$ on $\{ 0 \} \times \cdot \cdot \cdot \times \mathbb{R}^{n_j} \times \cdot \cdot \cdot \times \{ 0 \}$ for each $j=1,...,m$.}

\qquad

\textit{Proof.} It is easy to check that $$\bigoplus_{i=1}^{m} \bigcap_{1\leq j \neq i \leq m} \mathcal{N}_j = \mathbb{R}^{n} \Rightarrow
\bigoplus_{1\leq i \neq k \leq m} \bigcap_{1\leq j \neq i \leq m} \mathcal{N}_j = \mathcal{N}_k$$ so
\begin{equation}
A_{k} \left( \bigcap_{1\leq j \neq k \leq m} \mathcal{N}_j \right) = \mathcal{R}(A_k), \,\,\,\,\,\,\,\, k=1, ..., m, \label{proy}
\end{equation}
since $dim(\mathcal{N}_k) = n - n_k$ then $dim\left( \bigcap_{1\leq j \neq k \leq m} \mathcal{N}_j \right) = dim(\mathcal{R}(A_k)) = n_k$. Let $\{ \gamma_{1}^{k}, ..., \gamma_{n_{k}}^{k} \}$ be a basis of
$\bigcap_{1\leq j \neq k \leq m} \mathcal{N}_j$ thus $\{ \gamma_{1}^{1}, ..., \gamma_{n_{1}}^{1}, ..., \gamma_{1}^{m}, ..., \gamma_{n_{m}}^{m} \}$ is a basis for
$\mathbb{R}^{n}$. Let $C$ be the $n \times n$ matrix which columns are the elements of the above basis. Since $A_1+ ...+ A_m$ is invertible we have that
$B=(A_1+ ...+ A_m)C$ is invertible, so (\ref{proy}) gives that $B^{-1}A_j C$ is the canonical projection from $\mathbb{R}^{n}$ on $\{ 0 \} \times \cdot \cdot \cdot \times \mathbb{R}^{n_j} \times \cdot \cdot \cdot \times \{ 0 \}$ for each $j=1,...,m$. $\blacksquare$

\qquad

\section{The main result}

\qquad

\textit{Proof of the Theorem 1.} We begin by obtaining the  $L^{p}-L^{q}$ boundedness of the operator $T_r$ for $1 < p < \frac{1}{r}$ and
$\frac{1}{q}=\frac{1}{p}-r$, and then with this result we will prove the $H^{p}-L^{q}$ boundedness of $T_r$ for $0 < p \leq 1$ and  $\frac{1}{q}=\frac{1}{p}-r$.

\qquad

\textit{$L^{p}-L^{q}$ boundedness.} If $A$ is a $n \times n$ invertible matrix we put $f_{A}(x) = f(A^{-1}x)$. Let $B$ and $C$ be the matrices given by Lemma 2, then
$$\left[T_r \left( f_{C}\right) \right]_{B^{-1}} (x) = $$
$$\int_{\mathbb{R}^{n}} \, \left|( Bx - A_1 y ) \right|^{-n_1 + \alpha_1} \cdot \cdot \cdot \left|( Bx - A_m y ) \right|^{-n_m + \alpha_m} f(C^{-1}y) \, dy =$$
$$\left|det (C) \right| \int_{\mathbb{R}^{n}} \, \left|B (x - B^{-1}A_1 C y ) \right|^{-n_1 + \alpha_1} \cdot \cdot \cdot \left|B (x - B^{-1}A_m C y ) \right|^{-n_m + \alpha_m} f(y) \, dy.$$
Since $B$ is invertible, then there exists a positive constant $c$ such that $c|x| \leq |Bx|$ for all $x \in \mathbb{R}^{n}$. Thus
$$\left| \left[T_r \left( f_{C}\right) \right]_{B^{-1}} (x) \right| \leq c
\int_{\mathbb{R}^{n}} \, \left|x - B^{-1}A_1 C y \right|^{-n_1 + \alpha_1} \cdot \cdot \cdot \left|x - B^{-1}A_m C y \right|^{-n_m + \alpha_m} |f(y)| \, dy$$
$$\leq c \int_{\mathbb{R}^{n_1} \times ... \times \mathbb{R}^{n_m}} \, \left|x^{1} - y^{1} \right|^{-n_1 + \alpha_1} \cdot \cdot \cdot \left|x^{m} - y^{m} \right|^{-n_m + \alpha_m} |f(y^{1},...,y^{m})| \, dy^{1}...dy^{m},$$
the second inequality it follows from Lemma 2 and from that $|x^{j}-y^{j}| \leq |x - P_j y|$, where $P_j=B^{-1}A_j C$ is the canonical projection from $\mathbb{R}^{n}$ on $\{0\} \times \cdot \cdot \cdot \times \mathbb{R}^{n_j} \times \cdot \cdot \cdot \times \{0\}$. Since $\gamma(\alpha_j)^{-1}|x^{j}-y^{j}|^{-n_j+\alpha_j}$, for an appropriate constant $\gamma(\alpha_j)$ (see \cite{stein2}, p. 117), is the kernel of the Riesz potential on $\mathbb{R}^{n_j}$, then Theorem 1 (in \cite{stein2}, p. 119) and Lemma 1 give the $L^{p}-L^{q}$ boundedness of the operator $T_r$ for $1 < p < \frac{1}{r}$ and $\frac{1}{q}=\frac{1}{p}-r$.

\qquad

\textit{$H^{p}-L^{q}$ boundedness.} Let $0<p\leq 1.$ We recall that a $p$-atom is a measurable function $a$ supported on a ball $B$ of $\mathbb{R}^{n}$ satisfying

$a)$ $\left\Vert a\right\Vert _{\infty }\leq \left\vert B\right\vert^{-\frac{1}{p}}$

\qquad

$b)$ $\int y^{\beta }a(y)dy=0$ for every multiindex $\beta $ with $\left\vert \beta \right\vert \leq \lfloor n(p^{-1}-1) \rfloor$, ($\lfloor s \rfloor$ denotes the integer part of $s$).

\qquad

Let $0<p\leq 1  < p_0 < \frac{1}{r}$, $0 < r < 1$, and $\frac{1}{q} = \frac{1}{p} - r$. Given $f\in H^{p}(\mathbb{R}^{n}) \cap L^{p_0}(\mathbb{R}^{n})$, from Theorem 2, p. 107, in \cite{stein}  we have that there exist a sequence of nonnegative numbers $\{ \lambda_j \}_{j=1}^{\infty}$, a sequence of balls $B_{j}=B(z_j, \delta_j)$ centered at $z_j$ with radius $\delta_j$ and $p$ - atoms $a_j$ supported on $B_j$, satisfying
\begin{equation}
\sum\limits_{j\in \mathbb{N}}\left\vert \lambda _{j}\right\vert ^{p} \leq c\left\Vert f\right\Vert _{H^{p}}^{p}, \label{esti}
\end{equation}
such that $f$ can be decomposed as $f=\sum\limits_{j\in \mathbb{N}}\lambda _{j} a_{j},$ where the convergence is in $H^{p}$ and in $L^{p_0}$ (for the converge in $L^{p_0}$ see Theorem 5 in \cite{Rocha}). So the $H^{p}-L^{q}$ boundedness of $T_r$ will be proved if we show that there exists $c>0$ such that
\begin{equation} 
\left\Vert T_r a_j \right\Vert _{L^{q}}\leq c \label{uniform estimate}
\end{equation}
with $c$ independent of the $p$-atom $a_j$. Indeed,  since $f = \sum_{j=1}^{\infty} \lambda_j a_j$ in $L^{p_0}$ and $T_{r}$ is a $L^{p_0}-L^{\frac{p_0}{1- r p_0}}$ bounded operator, we have that $|T_{r} f (x)| \leq \sum_{j=1}^{\infty} \lambda_j |T_{r} a_j(x)|$ a.e.$x$,
this pointwise estimate, the inequality in (\ref{uniform estimate}), joint to the inequality
\[
\left( \sum_{j \in \mathbb{N}} | \lambda_j |^{\min\{1, q\}} \right)^{\frac{1}{\min\{1, q\}}} \leq \left( \sum_{j \in \mathbb{N}} | \lambda_j |^{p} \right)^{\frac{1}{p}}
\]
and (\ref{esti}) we obtain $\left\Vert T_r f\right\Vert_q \leq c \left\Vert f\right\Vert _{H^{p}}$, for all $f\in H^{p}(\mathbb{R}^{n}) \cap L^{p_0}(\mathbb{R}^{n})$, so the theorem follows from the density of $ H^{p}(\mathbb{R}^{n}) \cap L^{p_0}(\mathbb{R}^{n})$ in $H^{p}(\mathbb{R}^{n})$.

We will prove the estimate in (\ref{uniform estimate}). We define $D = \max_{1\leq i \leq m} \max_{|y|=1} |A_i (y)|$. Let $a_j$ be an $p$ - atom supported on a ball $B_j = B(z_j, \delta_j)$, for each $1 \leq i \leq m$ let $B_{ji}^{\ast} = B(A_i z_j, 4 D \delta_j)$. Since $T_r$ is bounded from $L^{p_{0}}\left( \mathbb{R}^{n}\right) $ into $L^{q_{0}}\left( \mathbb{R}^{n}\right) $ for $\frac{1}{q_{0}}=\frac{1}{p_{0}}-r,$ $1<p_{0}<\frac{1}{r},$ the H\"{o}lder inequality gives
\begin{equation*}
\int\limits_{\bigcup\limits_{1\leq i\leq m}B_{ji}^{\ast }}\left\vert
T_r a_j(x)\right\vert ^{q}dx\leq \sum\limits_{1\leq i\leq
m}\int\limits_{B_{ji}^{\ast }}\left\vert T_r a_j(x)\right\vert ^{q}dx
\end{equation*}
\begin{equation*}
\leq c\sum\limits_{1\leq i\leq m}\left\vert B_{ji}^{\ast }\right\vert ^{1- \frac{q}{q_{0}}}\left\Vert T_r a_j \right\Vert _{q_{0}}^{q}\leq c \delta_{j}^{n-\frac{
nq}{q_{0}}}\left\Vert a_j \right\Vert _{p_{0}}^{q}
\end{equation*}
\begin{equation*}
\leq c \delta_{j}^{n-\frac{nq}{q_{0}}}\left( \int\limits_{B_j}\left\vert
a_j \right\vert ^{p_{0}}\right) ^{\frac{q}{p_{0}}}\leq c \delta_{j}^{n-\frac{nq}{
q_{0}}} \delta_{j}^{-\frac{nq}{p}} \delta_{j}^{\frac{nq}{p_{0}}}=c.
\end{equation*}%

We denote $k(x,y)=\left\vert x-A_{1}y\right\vert ^{- n_1 + \alpha
_{1}}...\left\vert x-A_{m}y\right\vert ^{- n_m + \alpha _{m}}$ and we put $N-1 = \lfloor n(p^{-1}-1) \rfloor$. In view of the moment condition of $a_j$ we have
\begin{equation}
T_{r} a_j(x)=\int\limits_{B_j}k(x,y)a_j(y)dy=\int\limits_{B_j}\left( k(x,y)-q_{N, j}\left(
x,y\right) \right) a_j(y)dy,
\end{equation}%
\newline
where $q_{N, j}$ is the degree $N-1$ Taylor polynomial of the function $
y\rightarrow k(x,y)$ expanded around $z_j$. By the standard estimate of the
remainder term of the taylor expansion, there exists $\xi $ between $y$ and
$z_j$ such that
\[
\left\vert k(x,y)-q_{N, j}\left( x,y\right) \right\vert \lesssim \left\vert
y-z_j \right\vert ^{N}\sum\limits_{k_{1}+...+k_{n}=N}\left\vert \frac{%
\partial ^{N}}{\partial y_{1}^{k_{1}}...\partial y_{n}^{k_{n}}}k(x,\xi
)\right\vert
\]
\[
\lesssim \left\vert y-z_j \right\vert ^{N}\left(
\prod\limits_{i=1}^{m}\left\vert x-A_{i}\xi \right\vert ^{-n_i + \alpha
_{i}}\right) \left( \sum\limits_{l=1}^{m}\left\vert x-A_{l}\xi \right\vert
^{-1}\right) ^{N}.
\]
Now we decompose $\mathbb{R}^{n} = \bigcup_{i=1}^{m} B_{ji}^{\ast} \cup R_j$, where $R_j = \left( \bigcup_{i=1}^{m} B_{ji}^{\ast} \right)^{c}$, at the same time we decompose $R_j = \bigcup_{k=1}^{m} R_{jk}$ with
$$R_{jk} = \{ x \in R_{j} : |x - A_k z_j| \leq |x - A_i z_j| \,\, for \,\, all \,\, i \neq k \}.$$
If $x \in R_{j}$ then $|x - A_i z_j| \geq 4D \delta_j$, for all $i=1, ..., m$, since $\xi \in B_j$ it follows that $|A_i z_j - A_i \xi | \leq D \delta_j \leq \frac{1}{4} |x - A_i z_j|$ so
$$|x - A_i \xi| = |x - A_i z_j + A_i z_j - A_i \xi| \geq |x - A_i z_j| - |A_i z_j - A_i \xi| \geq \frac{3}{4} |x - A_i z_j|.$$
If $x \in R_j$, then $x \in R_{jk}$ for some $k$, and since $\sum_{i=1}^{m}(-n_i + \alpha_i) = - n(1-r)$ we obtain
$$
\left\vert k(x,y)-q_{N, j}\left( x,y\right) \right\vert  \lesssim
\left\vert y-z_j \right\vert ^{N}\left( \prod\limits_{i=1}^{m}\left\vert
x-A_{i}z_j \right\vert ^{- n_i + \alpha _{i}}\right) \left(
\sum\limits_{l=1}^{m}\left\vert x-A_{l}z_j \right\vert ^{-1}\right) ^{N}, 
$$
$$
\lesssim \left\vert y-z_j \right\vert ^{N} \left\vert x-A_{k}z_j \right\vert ^{-n(1-r) -N}, \,\,\,\,\, if \,\, x \in R_{jk} \,\,\, and \,\,\, y \in B_j.
$$
This inequality allow us to conclude that
\[
\int\limits_{R_j}\left\vert \int\limits_{B_j}K(x,y)a_j (y)dy\right\vert^{q}dx = \int\limits_{R_j}\left\vert \int\limits_{B_j} \left[ K(x,y) - q_{N,j}(x,y) \right] a_j (y)dy\right\vert^{q}dx 
\]
\[
\lesssim \sum_{k=1}^{m} \int\limits_{R_{jk}}\left(
\int\limits_{B_j}\left\vert y-z_j \right\vert ^{N} \left\vert x-A_{k}z_j \right\vert ^{-n (1-r) - N}\left\vert
a_j (y)\right\vert dy\right) ^{q}dx
\]
\[
\lesssim \left( \int\limits_{B_j}\left\vert y- z_j \right\vert ^{N}\left\vert
a_j (y)\right\vert dy\right) ^{q} \sum_{k=1}^{m} \int\limits_{\left(B_{jk}^{\ast} \right)^{c}}\left\vert x-A_{k}z_j \right\vert ^{-n (1-r)q - Nq} dx
\]
\[
\lesssim  \delta_{j}^{qN-n \frac{q}{p}+nq} \int\limits_{4 D \delta_j}^{\infty }t^{-q\left( n(1-r)+N\right) +n-1}dt \leq c
\]
with $c$ independent of the $p-$atom $a_j$ since $-q\left( n(1-r)+N\right)
+n<0.\blacksquare $

\section{An example}

For $n=m=3$, $n_1 = n_2 = n_3 = 1$ we consider the following $3 \times 3$ singular  matrices
\[ 
A_1= \left(
\begin{array}{ccc}
                      \,\,\, 4 & \,\,\, 4 & -1 \\
                      \,\,\, 0 & \,\,\,0 & \,\,\, 0 \\
                      -4 & -4 & \,\,\, 1 \\
                    \end{array}
                  \right), \,\,\,\,\,\, A_2=\left(
\begin{array}{ccc}
                      \,\,\, 1 & -1 & 0 \\
                      -2 & \,\,\, 2 & 0 \\
                      \,\,\, 0 & \,\,\, 0 & 0 \\
                    \end{array}
                  \right), \,\,\,\,\,\, A_3 = \left(
                                      \begin{array}{ccc}
                                        \,\,\, 1 & 0 & -1 \\
                                        -3 & 0 & \,\,\, 3 \\
                                        -1 & 0 & \,\,\, 1 \\
                                      \end{array}
                                    \right).
\]

It is clear that 
\[
A_1 + A_2 +A_3 =  \left(
                                      \begin{array}{ccc}
                                        \,\,\, 6 & \,\,\, 3 & -2 \\
                                        -5 & \,\,\, 2 & \,\,\, 3 \\
                                        -5 & -4 & \,\,\, 2 \\
                                      \end{array}
                                    \right)
\]
is invertible. For each $j=1, 2, 3,$  let $\mathcal{N}_j = \{ x \in \mathbb{R}^{3} : A_j x = 0\}$. A computation gives
$$
\mathcal{N}_1 = \langle (1, 0, 4) , (0, 1, 4) \rangle, \,\,\, \mathcal{N}_2 = \langle (1, 1, 0), (0, 0, 1) \rangle, \,\,\, \mathcal{N}_3 = \langle (1, 0, 1), (0, 1, 0) \rangle,
$$
one can check that
$$
\mathcal{N}_1  \cap \mathcal{N}_2 = \langle (1, 1, 8) \rangle, \,\,\, \mathcal{N}_1  \cap \mathcal{N}_3 = \langle (4, -3, 4) \rangle,
\,\,\, and \,\,\, \mathcal{N}_2  \cap \mathcal{N}_3 = \langle (1, 1, 1) \rangle.
$$

\qquad

We observe that $\mathcal{N}_1 \cap \mathcal{N}_2  \oplus \mathcal{N}_1  \cap \mathcal{N}_3 \oplus  \mathcal{N}_2  \cap \mathcal{N}_3 = \mathbb{R}^{3}$. As in the proof of Lemma 2, we define the matrices $C$ and $B$ by
\[
C = \left(
                \begin{array}{ccc}
                 1 & \,\,\, 4 & 1 \\
                 1 & -3 & 1 \\
                 1 & \,\,\, 4 & 8\\
                  \end{array}
                  \right),
\,\,\, B= ( A_1 + A_2 + A_3 ) \,  C = \left(
                                      \begin{array}{ccc}
                                        \,\,\, 7 & \,\,\, 7 & -7 \\
                                        \,\,\, 0 & -14 & \,\,\, 21 \\
                                        -7 & \,\,\, 0 & \,\,\, 7 \\
                                      \end{array}
                                    \right), 
\]
both matrices are invertibles with
\[
B^{-1} = \left(
                                      \begin{array}{ccc}
                                        \frac{2}{21} & \frac{1}{21} & -\frac{1}{21} \\\\
                                        \frac{1}{7} & 0 & \,\,\, \frac{1}{7} \\\\
                                        \frac{2}{21} & \frac{1}{21} & \,\,\, \frac{2}{21} \\
                                      \end{array}
                                    \right).
\]
Now it is easy to check that
\[
B^{-1 } A_1 C = \left(
                                      \begin{array}{ccc}
                                        1 & 0 & 0 \\
                                        0 & 0 & 0 \\
                                        0 & 0 & 0 \\
                                      \end{array}
                                    \right), \,\,\, 
B^{-1 } A_2 C = \left(
                                      \begin{array}{ccc}
                                        0 & 0 & 0 \\
                                        0 & 1 & 0 \\
                                        0 & 0 & 0 \\
                                      \end{array}
                                    \right), \,\,\, 
B^{-1 } A_3 C = \left(
                                      \begin{array}{ccc}
                                        0 & 0 & 0 \\
                                        0 & 0 & 0 \\
                                        0 & 0 & 1 \\
                                      \end{array}
                                    \right). 
\]
So, from Theorem 1, it follows that the operator $T_r$ defined by
\[
T_rf(x) = \int_{\mathbb{R}^{3}} |x-A_1 y |^{-1+r} |x -A_2 y |^{-1+r} |x - A_3 y |^{-1+r} f(y) dy,
\]
 with $0 < r < 1$, is a bounded operator from $H^{p}(\mathbb{R}^{3})$ into $L^{q}(\mathbb{R}^{3})$ for $0 < p < 1/r$ and $\frac{1}{q} = \frac{1}{p} - r$.
\qquad

\end{document}